\documentclass[12pt]{amsart}    
\usepackage{amsmath, amssymb} 
\usepackage{graphicx} 

\newtheorem{theorem}{Theorem}[section]
\newtheorem{proposition}[theorem]{Proposition}
\newtheorem{lemma}[theorem]{Lemma}

\theoremstyle{definition}

\newtheorem{remark}[theorem]{Remark}

\numberwithin{equation}{section}

\begin{document}

%

\title[]{On rectangular diagrams, \\ 
Legendrian knots and transverse knots} 

\author[]{Hiroshi MATSUDA} 

\author[]{William W. MENASCO}

\address{Department of Mathematics, 
Graduate School of Science, 
Hiroshima University, 
Hiroshima 739-8526, JAPAN} 
\curraddr{Department of Mathematics, 
Columbia University, 
New York, NY 10027, USA} 

\email{matsuda@math.sci.hiroshima-u.ac.jp, matsuda@math.columbia.edu} 

\address{Department of Mathematics, 
University at Buffalo, Buffalo, NY 14260} 

\email{menasco@math.buffalo.edu}

\maketitle

 
\section{introduction} 
 
A correspondence is studied in \cite{matsuda}
between front projections of 
Legendrian links in $({\mathbb{R}}^3, \xi_{\rm std})$ 
and rectangular diagrams. 
In this paper, 
we introduce braided rectangular diagrams, and 
study a relationship with Legendrian links 
in $({\mathbb{R}}^3, \xi_{\rm sym})$. 
We show Alexander and Markov Theorems 
for Legendrian links in $({\mathbb{R}}^3, \xi_{\rm sym})$.

\begin{figure}
 \includegraphics[keepaspectratio]{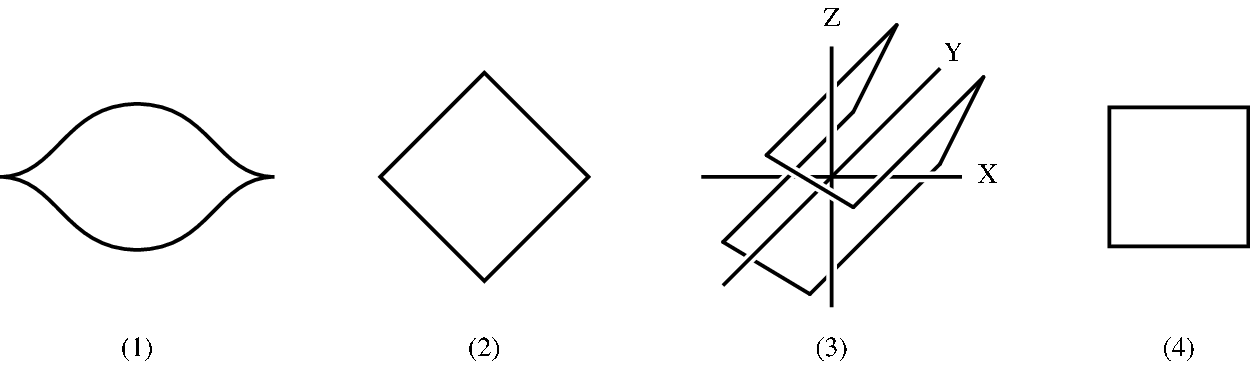}
\caption{} 
\label{trivial-knot}
\end{figure}

We review a relationship between front projections of 
Legendrian links in $({\mathbb{R}}^3, \xi_{\rm std})$ 
and rectangular diagrams. 
The standard contact structure $\xi_{\rm std}$ in ${\mathbb{R}}^3$ 
is defined by $\ker (dz - y dx)$. 
It is well-known \cite{swiatkowski} that 
every Legendrian link in $({\mathbb{R}}^3, \xi_{\rm std})$ has 
a front projection with transverse double points and 
cusp singularities. 
Figure \ref{trivial-knot} (1) illustrates a front projection 
of a topologically trivial Legendrian knot. 
Changing every point with a horizontal tangent 
in a front projection to a corner, Figure \ref{trivial-knot} (2), 
followed by rotating the obtained diagram 45 degree clockwise, 
we obtain a rectangular diagram, Figure \ref{trivial-knot} (4).

The intersection of the plane $\{ (x, y, z) \ | \ y = a \}$ 
with contact planes, called a {\it characteristic foliation}, 
consists of lines $z = ax + {\rm (constant)}$, 
where $a \in {\mathbb{R}}$. 
The characteristic foliation on the plane $\{ (x, y, z) \ | \ x = b \}$ 
consists of lines $z = {\rm (constant)}$, where $b \in {\mathbb{R}}$. 
Changing every vertical arc in a rectangular diagram 
to a subarc of the line $\{ (x, -1, z) \ | \ z = -x + v_i \}$ 
for some $v_i \in {\mathbb{R}}$, 
every horizontal arc to a subarc of the line 
$\{ (x, 1, z) \ | \ z = x + h_j \}$ for some $h_j \in {\mathbb{R}}$, and 
every corner to a subarc of the line 
$\{ (x_i, y, z_j) \ | \ y \in {\mathbb{R}}\}$ 
for some $x_i \in {\mathbb{R}}$ and $z_j \in {\mathbb{R}}$, 
we obtain a Legendrian link in $({\mathbb{R}}^3, \xi_{\rm std})$ 
from a rectangular diagram. 
We notice that the Legendrian link 
constructed as above from a rectangular diagram 
is Legendrian isotopic to the Legendrian link 
corresponding to the front projection. 
For example, 
the diagram in Figure \ref{trivial-knot} (4) corresponds to 
a Legendrian knot, illustrated in  Figure \ref{trivial-knot} (3),  
consisting of the following eight arcs: \\ 
$\{ (x, -1, z) \ | \ x+z=1, -\varepsilon \leq x \leq 1-\varepsilon, 
\varepsilon \leq z \leq 1+\varepsilon \}$, 
$\{ (-\varepsilon, y, 1+\varepsilon) \ | \ -1 \leq y \leq 1 \}$, \\ 
$\{ (x, 1, z) \ | \ x-z=-1, -1-\varepsilon \leq x \leq -\varepsilon, 
\varepsilon \leq z \leq 1+\varepsilon \}$, 
$\{ (-1-\varepsilon, y, \varepsilon) \ | \ -1 \leq y \leq 1 \}$, \\ 
$\{ (x, -1, z) \ | \ x+z=-1, -1-\varepsilon \leq x \leq -\varepsilon, 
-1+\varepsilon \leq z \leq \varepsilon \}$, 
$\{ (-\varepsilon, y, -1+\varepsilon) \ | \ -1 \leq y \leq 1 \}$, \\ 
$\{ (x, 1, z) \ | \ x-z=1, -\varepsilon \leq x \leq 1-\varepsilon, 
-1+\varepsilon \leq z \leq \varepsilon \}$, 
$\{ (1-\varepsilon, y, \varepsilon) \ | \ -1 \leq y \leq 1 \}$, \\ 
where $\varepsilon$ is a small positive number.
The union of these arcs is piecewise Legendrian.  We obtain a Legendrian knot by Legendrian-smoothing
our edgepath in arbitrarily small neighborhoods around the endpoints of each arc.

In $\S$2, we study a similar correspondence 
between braided rectangular diagrams and 
Legendrian links in $({\mathbb{R}}^3, \xi_{\rm sym})$. 
This leads us to Alexander and Markov Theorems 
for Legendrian links in $({\mathbb{R}}^3, \xi_{\rm sym})$. 
This answers Problem 2 in \cite{orevkov-shevchishin}, 
implicitly stated also in \cite{epstein-fuchs-meyer}. 
Alexander and Markov Theorems for transverse links 
in $({\mathbb{R}}^3, \xi_{\rm sym})$ was proved 
in \cite{bennequin}, 
\cite{orevkov-shevchishin} and \cite{wrinkle}.

\begin{theorem} 
\label{transverse-alexander} 
\cite{bennequin}
{\rm (Alexander Theorem for transverse links 
in $({\mathbb{R}}^3, \xi_{\rm sym})$)} \\ 
Any transverse link in $({\mathbb{R}}^3, \xi_{\rm sym})$ 
is transversely isotopic to a closed braid. 
\end{theorem}

\begin{theorem} 
\label{transverse-markov} 
\cite{orevkov-shevchishin}, \cite{wrinkle} 
{\rm (Markov Theorem for transverse links 
in $({\mathbb{R}}^3, \xi_{\rm sym})$)}  \\ 
Two closed braids represent the same transverse link 
if and only if they are related by positive stabilizations and 
conjugation in the braid group. 
\end{theorem}

In $\S$3, 
we describe a construction of 
an explicit example of the Etnyre-Honda pair of 
Legendrian knots \cite{etnyre-honda}, 
which is announced in \cite{menasco-addendum}.

Acknowledgement.  
The first author is partially supported by 
JSPS Postdoctoral Fellowship for Research Abroad.
This work started 
when first author visited SUNY at Buffalo in September-October 2003, 
and was completed 
the second author's visited Hiroshima University in June 2005.
The second author was partially supported by NSF grant \#DMS 0306062.
Parts of this work were presented 
by the first author in April 2005 at Hiroshima University.

\section{Alexander and Markov Theorems for 
Legendrian links in $({\mathbb{R}}^3, \xi_{\rm sym})$} 

Let $\{ H_\theta \ | \ 0 \leq \theta < 2\pi \}$ denote 
an open book decomposition of ${\mathbb{R}}^3$, 
that is, ${\mathbb{R}}^3 \setminus \{ z$-axis$\}$ is fibered 
by a collection of half-plane fibers $H_\theta$, 
where the boundary of $H_\theta$ is the $z$-axis. 
When we use the cylindrical coordinates $(r, \theta, z)$ 
on ${\mathbb{R}}^3$, 
a fiber $H_{\theta_0}$ is the set 
$\{ (r, \theta, z) \ | \ \theta = \theta_0 \}$. 
An oriented link $X$ in ${\mathbb{R}}^3$ is 
a {\it closed} $n$-{\it braid} 
if $X \subset {\mathbb{R}}^3 \setminus \{ z$-axis$\}$ 
intersects each fiber $H_\theta$ transversely and 
positively in $n$ points. 
Possibly after a small isotopy of $X$ in 
${\mathbb{R}}^3 \setminus \{ z$-axis$\}$, 
we can consider a regular projection 
$\pi \colon X \to C_1$ given by 
$(r, \theta, z) \mapsto (1, \theta, z)$, 
where $C_1 = \{ (r, \theta, z) \ | \ r = 1 \}$. 
We may assume that the singularities consist of $\pi(X)$ consists of 
finitely many transverse double points.

A {\it horizontal arc}, $h \subset C_1$, is an arc 
with a parameterization 
$\{ (1, \theta(t), z_0) \ | \ 
0 \leq t \leq 1, \theta(t) \in [\theta_1, \theta_2] \}$, 
where $| \theta_1 - \theta_2 | < 2\pi$. 
The {\it horizontal level} of $h$ is a fixed constant $z_0$, and 
the {\it angular support} of $h$ is the interval $[\theta_1, \theta_2]$. 
Horizontal arcs inherit a natural orientation from the forward 
direction of the $\theta$ coordinate. 
A {\it vertical arc}, $v \subset H_{\theta_0}$, 
is an arc with a parameterization $\{ (r(t), \theta_0, z(t)) \ | \ 
0 \leq t \leq 1, r(0) = r(1) = 1;$ and $r(t) > 1, \frac{dz(t)}{dt} \neq 0$ 
for $t \in (0, 1) \}$, 
where $r(t)$ and $z(t)$ are ${\mathbb{R}}$-valued functions 
that are continuous on $[0, 1]$ and differential on $(0, 1)$. 
The {\it angular level} of $v$ is $\theta_0$, 
the {\it vertical support} of $v$ is 
the interval $[z(0), z(1)]$ or $[z(1), z(0)]$.

An oriented link $X \subset {\mathbb{R}}^3$ is 
an {\it arc presentation} \cite{menasco-exchange} 
if $X = h_1 \cup v_1 \cup \cdots \cup h_n \cup v_n$ 
satisfies the following conditions: \\ 
(1) each $h_i$, $(i = 1, \cdots, n)$, is an oriented horizontal arc 
with its inherited orientation agreeing with the orientation of $X$, \\ 
(2) each $v_i$, $(i = 1, \cdots, n)$, is a vertical arc, \\ 
(3) the intersection $h_i \cap v_j$ consists of 
a point of $\partial h_i \cap \partial v_j$ for 
$i = 1, \cdots, n$ and $j \equiv \{ i - 1, i \} \pmod{n}$, 
and $h_i \cap v_j = \emptyset$ if $j \not\equiv \{ i - 1, i \} \pmod{n}$, \\ 
(4) the horizontal level of each horizontal arc is distinct, and 
the angular level of each vertical arc is distinct, \\ 
(5) the orientations of vertical arcs are assigned 
so as to make the components of $X$ oriented.

A projection of an arc presentation of $X$ onto $C_1$ 
with over/under informations is called 
a {\it braided rectangular diagram}. 
A projection $\pi(X)$ without the conditions about orientations 
in (1) and (5) is called a {\it rectangular diagram}. 
It is proved in \cite{dynnikov} and \cite{menasco-exchange} that 
every link in ${\mathbb{R}}^3$ has a braided rectangular diagram. 

By a contactmorphism $({\mathbb{R}}^3, \xi_{\rm std})$ is equivalent to the symmetric
contact structure $({\mathbb{R}}^3, \xi_{\rm sym})$ where $\xi_{\rm sym} = {\rm ker} \alpha_{\rm sym}$
for $\alpha_{\rm sym} = dz + x dy - ydx$ (in Euclidean coordinates) or $\alpha_{sym} = dz + r^2 d \theta$
(in clyindrical coordinates).  The symmetric contact structure was the structure utilized
by D. Bennequin \cite{bennequin} in his classical argument the any transversal knot is transversally
isotopic to a braid---the $z$-axis being the designated braid axis.

A rectangular diagram on $C_1$ corresponds to 
a Legendrian link in $({\mathbb{R}}^3, \xi_{\rm sym})$ as follows. 
First we notice that 
the characteristic foliation on 
the cylinder $\{ (r, \theta, z) \ | \ r = r_0 \}$ 
consists of $``$spiral$"$ curves of slope $\frac{dz}{d\theta} = - r_0^2$, 
and that the characteristic foliation on the plane 
$\{ (r, \theta, z) \ | \ \theta = \theta_0 \}$ 
consists of lines $\{ (r, \theta_0, z_0) \ | \ r > 0 \}$ 
for $z_0 \in {\mathbb{R}}$. 
Change every horizontal arc in a rectangular diagram 
to a subarc, a {\it near-horizontal arc}, of 
the characteristic foliation on the cylinder 
$\{ (r, \theta, z) \ | \ r = r_1 \}$ with $r_1$ sufficiently small, 
and every vertical arc to a subarc, a {\it near-vertical arc}, 
of the characteristic foliation 
on the cylinder $\{ (r, \theta, z) \ | \ r = r_2 \}$ 
with $r_2$ sufficiently large. 
Change every corner in a rectangular diagram 
to a subarc of the line $\{ (r, \theta, z) \ | \ \theta = \theta_0, r > 0 \}$ 
which connects one endpoint of a near-horizontal arc on the cylinder 
$\{ (r, \theta, z) \ | \ r = r_1 \}$ and one endpoint of a near-vertical arc 
on the cylinder $\{ (r, \theta, z) \ | \ r = r_2 \}$. 
Adjusting the $r$-coordinates of near-horizontal and 
near-vertical arcs properly and Legendrian smoothing in arbitrarily small neighborhoods of the
arc endpoints, 
we obtain a Legendrian link $L$ in $({\mathbb{R}}^3, \xi_{\rm sym})$ 
from a braided rectangular diagram. 
The vertical (resp. angular) support of a near-horizontal 
(resp. near-vertical) arc of $L$ is the interval 
in the $z$-coordinate (resp. $\theta$-coordinate) containing the arc. 
We may  Legendrian isotope $L$ so that: \\ 
(1) the vertical support of a near-horizontal arc of $L$ 
is sufficiently small, and is disjoint from each other, \\ 
(2) the angular support of a near-vertical arc of $L$ 
is sufficiently small, and is disjoint from each other. 

A Legendrian link $L$ has a {\it horizontal/vertical disjoint property} 
if the above conditions are satisfied. 

Next we construct a rectangular diagram 
from a Legendrian link in $({\mathbb{R}}^3, \xi_{\rm sym})$. 

\begin{lemma} 
Let $L$ be a Legendrian link in $({\mathbb{R}}^3, \xi_{\rm sym})$. 
Then $L$ may be Legendrian isotoped to $L^\prime$ so that 
$L^\prime$ in in the half-space with $y > 0$.  In particular, $L^\prime \cap H_\pi = \emptyset$. 
\end{lemma}

\begin{proof} 
Instead of Legendrian isotoping $L$, 
we describe a contactomorphism of 
$({\mathbb{R}}^3, \xi_{\rm sym})$ that takes $L$ to $L'$. 
Let $f$ be a diffeomorphism of ${\mathbb{R}}^3$ 
defined by $(x, y, z) \mapsto (x+ {\rm K}, y+ {\rm K}, 
z+ {\rm K} (x-y))$, where ${\rm K} \in {\mathbb{R}}$.
Starting with the Euclidean version of $\alpha_{\rm sym}$, we then
obtain $({\mathbb{R}}^3, \xi_{\rm sym}^\prime)$, 
where $\xi_{\rm sym} '$ is defined by the kernel of the 1-form 
$\alpha_{\rm sym}^\prime = d(z+{\rm K}(x-y)) +  (x + {\rm K}) d(y+ {\rm K}) - (y+{\rm K}) d(x+{\rm K}) 
= dz + x dy - y dx$. 
By the compactness of $L$, for a large enough ${\rm K}$
the image of $L$ by $f$, $L^\prime$, will be in the half-space $y>0$.  In particular,
$L^\prime \cap H_\pi = \emptyset$.
\end{proof}

Let $L$ be a Legendrian link in $({\mathbb{R}}^3, \xi_{\rm sym})$ 
that is in the half-space having $y>0$.  Then it does not intersect the $z$-axis and 
the image of $L$ onto $C_1$ by 
$\pi \colon (r, \theta, z) \mapsto (1, \theta, z)$, $\pi(L)$ is well-defined.
We called $\pi(L)$ the {\it cylinder projection} of $L$. 
We notice that when $L \cap H_\pi = \emptyset$ the front projections of Legendrian links 
in $({\mathbb{R}}^3, \xi_{\rm std})$ and
cylinder projections of Legendrian links in 
$({\mathbb{R}}^3, \xi_{\rm sym})$ are combinatorial equivalent.
We can then adapt the proof of the Reidemeister Theorem for Legendrian links 
in terms of front projections to the setting of cylinder projections to show the following. 
See the proof of Theorem B in \cite{swiatkowski}.

\begin{figure} 
 \includegraphics[keepaspectratio]{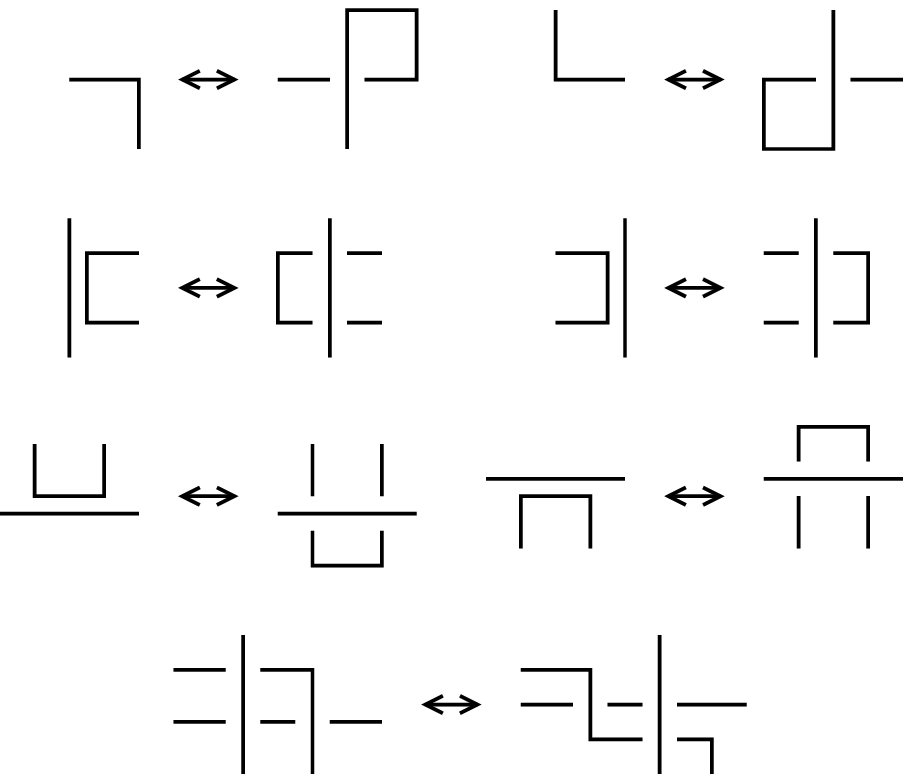}
\caption{} 
\label{reidemeister}
\end{figure}

\begin{proposition} 
Let $L_1$ and $L_2$ be Legendrian links in 
$({\mathbb{R}}^3, \xi_{\rm sym})$ such that 
each of $L_1$ and $L_2$ are in the half-space having $y>0$. 
Let $D_1$ and $D_2$ be the cylinder projections of 
$L_1$ and $L_2$, respectively. 
Legendrian links $L_1$ and $L_2$ are Legendrian isotopic 
in ${\mathbb{R}}^3 \setminus \{ z-axis \}$ if and only if 
$D_1$ and $D_2$ are related by regular homotopy and 
a finite sequence of moves 
that are obtained from the diagrams in Figure \ref{reidemeister} 
by rotating $\Theta$ degree counterclockwise with $0 < \Theta < 90$. 
\end{proposition}

As we obtain a rectangular diagram from a front projection, 
we obtain a rectangular diagram on $C_1$ from a cylinder projection 
by changing every point with $\frac{dz}{d\theta} = -1$ to a corner.

\begin{figure} 
 \includegraphics[keepaspectratio]{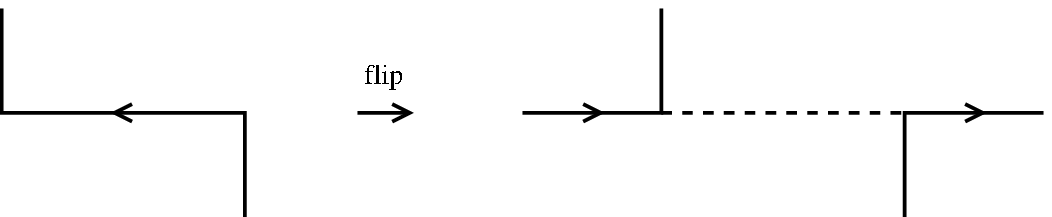}
\caption{} 
\label{flip}
\end{figure}

Next we define an operation, called a {\it flip}.  
Let $L = \cup_{i=1}^n (h_i \cup c(h_i) \cup v_i \cup c(v_i))$ be 
a Legendrian link that corresponds to a rectangular diagram, 
where $h_i$ (resp. $v_i$) is 
a near-horizontal (resp. near-vertical) arc, 
and $c(h_i)$ connects $h_i$ and $v_i$, and 
$c(v_i)$ connects $v_i$ and $h_{i+1}$. 
We may assume that 
a Legendrian link $L$ has a horizontal/vertical disjoint property.
Moreover, by a Legendrian isotopy that corresponds to scaling the
the angle between $v_{i-1}$ and $v_i$ we can assume that the angle between
the angular support of $h_i$ is $\pi$.
We may Legendrian isotope $L$ so that 
the $r$-coordinate of $h_i$ goes to 0, and that 
the $\theta$-coordinates of $c(v_{i-1})$ and $c(h_i)$ remain fixed. 
Then we have to adjust the lengths of $c(v_{i-1})$ and $c(h_i)$, 
and the $r$-coordinates of $v_{i-1}$ and $v_i$, and 
the lengths of $c(h_{i-1})$ and $c(v_i)$. 
After these Legendrian isotopies of $L$, 
we may assume that $L$ has a horizontal/vertical disjoint property. 
We may Legendrian isotope $L$ so that 
$h_i$ shrinks to a point on the $z$-axis. 
Therefore the $z$-coordinate of $c(v_{i-1})$ and $c(h_i)$ 
are the same. 
We may further Legendrian isotope $L$ to $L'$ 
so that $h_i$ passes through the $z$-axis to $h_i '$, 
that (the angular support of $h_i$) $\cap$ 
(the angular support of $h_i '$) $=$ 
(the $\theta$-coordinates of $c(v_{i-1})$ and $c(h_i)$), 
where $h_i '$ is a subarc of $L'$.  By our assumption that the angular support of $h_i$ is
$\pi$ this isotopy wil be Legendrian as $L$ passes through the $z$-axis.
This isotopy from $L$ to $L'$ is called a {\it flip}, and 
the above argument shows that 
a flip is a Legendrian isotopy from $L$ to $L'$.

\begin{theorem} 
\label{legendrian-alexander} 
{\rm (Alexander Theorem for Legendrian links 
in $({\mathbb{R}}^3, \xi_{\rm sym})$)} \\ 
Every Legendrian link in $({\mathbb{R}}^3, \xi_{\rm sym})$ 
is Legendrian isotopic to a Legendrian link constructed from 
a braided rectangular diagram. 
\end{theorem}

\begin{proof} 
Let $D$ be a rectangular diagram of a Legendrian link $L$ 
with a horizontal/vertical disjoint property. 
Let $h$ be a near-horizontal arc of $L$ 
such that the induced orientation of $h$ from that of $L$ 
disagrees with the forward direction of the $\theta$-coordinate. 
We apply a flip to every such near-horizontal arc $h$. 
Then $L$ corresponds to a braided rectangular diagram. 
\end{proof}

When a Legendrian link $L$ intersects the $z$-axis 
during Legendrian isotopy of $L$ 
in $({\mathbb{R}}^3, \xi_{\rm sym})$, 
a neighborhood of the intersection of $L$ with the $z$-axis 
is horizontal. 
Therefore we may assume that 
their rectangular diagrams are related by one flip, 
so we have the following.

\begin{theorem} 
\label{legendrian-reidemeister} 
{\rm (Reidemeister Theorem for Legendrian links 
in terms of rectangular diagrams)} \\ 
Let $L_1$ and $L_2$ be Legendrian links in 
$({\mathbb{R}}^3, \xi_{\rm sym})$, and 
let $D_1$ and $D_2$ be the rectangular diagrams of 
$L_1$ and $L_2$, respectively. 
Legendrian links $L_1$ and $L_2$ are Legendrian isotopic 
in $({\mathbb{R}}^3, \xi_{\rm sym})$ if and only if 
$D_1$ and $D_2$ are related by 
a finite sequence of flips and 
moves in Figure \ref{reidemeister} on $C_1$. 
\end{theorem}

Next we see how positive and negative transverse push-offs 
are obtained from a Legendrian link corresponding to 
a braided rectangular diagram. 
Let $L=\cup_{i=1}^n (h_i \cup c(h_i) \cup v_i \cup c(v_i))$ 
be a Legendrian link in $({\mathbb{R}}^3, \xi_{\rm sym})$, 
where $h_i$ (resp. $v_i$) denotes 
a near-horizontal (resp. near-vertical) arc. 
Let $\alpha$ be a subarc of $L$ that is contained in a cylinder 
$C_{r_1} = \{ (r, \theta, z) \ | \ r = r_1 \}$,  
so $\alpha$ is $h_i$ or $v_i$. 
Let $[\theta_1, \theta_2]$ be the angular support of $\alpha$. 
Let $z(\theta)$ be the $z$-coordinate at $\theta$ of the line 
in the characteristic foliation on $C_{r_1}$ containing $\alpha$, 
defined in the $\theta$-interval 
$[\theta_1 - \varepsilon, \theta_2 + \varepsilon]$, 
where $\varepsilon$ is a small positive number.  
Let $\delta(\alpha)$ be a neighborhood of $\alpha$ on $C_{r_1}$ 
described as 
$\delta(\alpha) = \{ (r_1, \theta, z) \ | \ 
\theta_1- \varepsilon \leq \theta \leq \theta_2 + \varepsilon, 
z(\theta) -  \varepsilon \leq z \leq z(\theta) + \varepsilon \}$.  
Let $\Delta(\alpha)$ be a neighborhood of $\alpha$ 
in ${\mathbb{R}}^3$ described as 
$\Delta(\alpha) = 
[r_1 - \varepsilon, r_1 + \varepsilon] \times \delta(\alpha) = 
\{ (r, \theta, z) \ | \ 
r_1 - \varepsilon \leq r \leq r_1 + \varepsilon, 
\theta_1 - \varepsilon \leq \theta \leq \theta_2 + \varepsilon, 
z(\theta) - \varepsilon \leq z \leq z(\theta) + \varepsilon \}$.

Let $\beta$ be a subarc of $L$ which is contained in a line 
$\ell(\theta_1, z_1) = \{ (r, \theta, z) \ | \ 
\theta = \theta_1, z = z_1 \}$, that is, 
$\beta = \{ (r, \theta_1, z_1) \ | \ r_1 \leq r \leq r_2 \}$, 
so $\beta$ is $c(h_i)$ or $c(v_i)$. 
Let $\beta_\varepsilon =\{ (r, \theta_1, z_1) \ | \ 
r_1 - \varepsilon \leq r \leq r_2 + \varepsilon \}$ be 
a neighborhood of $\beta$ in $\ell(\theta_1, z_1)$. 
Let $\delta(\beta_\varepsilon, r)$ be 
a neighborhood of $\beta_\varepsilon \cap C_r$ on $C_r$ 
described as $\delta(\beta_\varepsilon, r_0)= 
\{ (r_0, \theta, z) \ | \ 
\theta_1 - \varepsilon \leq \theta \leq \theta_1 + \varepsilon, 
z(\theta) - \varepsilon \leq z \leq z(\theta) + \varepsilon \}$, 
where $z(\theta)$ denotes the $z$-coordinate of the integral line 
in the characteristic foliation on $C_{r_0}$ 
containing the point $\beta_\varepsilon \cap C_{r_0}$, 
defined in the $\theta$-interval 
$[\theta_1 - \varepsilon, \theta_2 + \varepsilon]$. 
Let $\Delta(\beta)$ be a neighborhood of $\beta$ in ${\mathbb{R}}^3$ 
described as 
$\Delta(\beta) = 
[r_1 - \varepsilon, r_2 + \varepsilon] \times \delta(\beta_\varepsilon, r) 
= \{ (r, \theta, z) \ | \ 
r_1 - \varepsilon \leq r \leq r_2 + \varepsilon, 
\theta_1 - \varepsilon \leq \theta \leq \theta_1 + \varepsilon, 
z(\theta) - \varepsilon \leq z \leq z(\theta) + \varepsilon \}$.

Let $\Delta(L)$ be a neighborhood of $L$ constructed
as $\Delta(L) = \cup_{i=1}^{2n} (\Delta(\alpha_i) \cup \Delta(\beta_i))$, 
where $L= \cup_{i=1}^{2n} (\alpha_i \cup \beta_i)$. 
After a small isotopy of $\Delta(L)$, 
we may assume that 
$\partial \Delta(L)$ consists of four sides. 
Two of the four sides of $\partial \Delta (L)$ 
each contain one Legendrian divide. 
Let $T_+(L)$ and $T_-(L)$ be the center lines of the other two sides. 
Then $T_+(L)$ and $T_-(L)$ are 
transverse push-offs of $L$ in $({\mathbb{R}}^3, \xi_{\rm sym})$. 
We note that $T_+(L)$ (resp. $T_-(L)$) is 
a positive (resp. negative) transverse push-off of $L$.

\begin{figure} 
 \includegraphics[keepaspectratio]{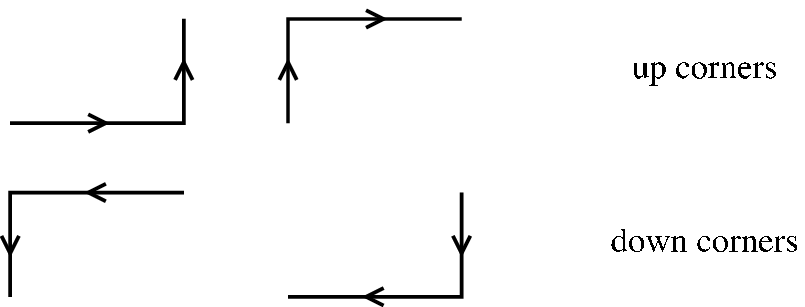}
\caption{} 
\label{corner}
\end{figure}

\begin{lemma} 
Let $L$ be a Legendrian link in $({\mathbb{R}}^3, \xi_{\rm sym})$ 
represented as a rectangular diagram $D$. 
Then we have $tb(L) = \omega(D) - \frac{1}{2} (d(D) + u(D))$ and 
$r(L) = n(D) + \frac{1}{2} (d(D) - u(D))$, where 
$u(D)$ (resp. $d(D)$) denotes 
the number of up (resp. down) corners of $D$, 
illustrated in Figure \ref{corner}, and 
$n(D)$ is the algebraic winding number of $L$ around the $z$-axis, 
and $\omega(D)$ is the algebraic crossing number of $D$ on $C_1$. 
\end{lemma}

\begin{proof} 
Let $v = \frac{\partial}{\partial z}$ be 
a vector field on ${\mathbb{R}}^3$. 
For any Legendrian knot $L$, 
$v$ is a vector field transverse to $\xi_{\rm sym}$ along $L$. 
So we have $tb(L) = \omega(D) - \frac{1}{2} (d(D) + u(D))$.

Let $w = \frac{\partial}{\partial r}$ be a vector field 
on $\{ x \geq \varepsilon \}$, 
$w = - \frac{\partial}{\partial r}$ a vector field 
on $\{ x \leq - \varepsilon \}$. 
We define a vector field $w$ on $\{ - \varepsilon < x < \varepsilon \}$ 
by interpolating between these two choices by rotating clockwise 
in the contact planes. 
We may Legendrian isotope $L$ 
so that all the vertical arcs of $D$ are contained 
in $\{ x \geq \varepsilon \}$. 
Then we have $r (L) = \frac{1}{2} (2n(D) + d(D) - u(D)) 
= n(D) + \frac{1}{2}(d(D) - u(D))$. 
\end{proof}

Similar arguments as above prove the following.

\begin{lemma} 
Let $T_+(L)$ be a positive transverse push-off 
of a Legendrian knot $L$ in $({\mathbb{R}}^3, \xi_{\rm sym})$ 
corresponding to a braided rectangular diagram $D$. 
Then we have $s\ell(T_+(L)) = \omega(D) - n(D)$. 
\end{lemma} 

Proposition 4 in \cite{dynnikov} and 
Theorem \ref{legendrian-reidemeister} 
prove the following.

\begin{figure} 
 \includegraphics[keepaspectratio]{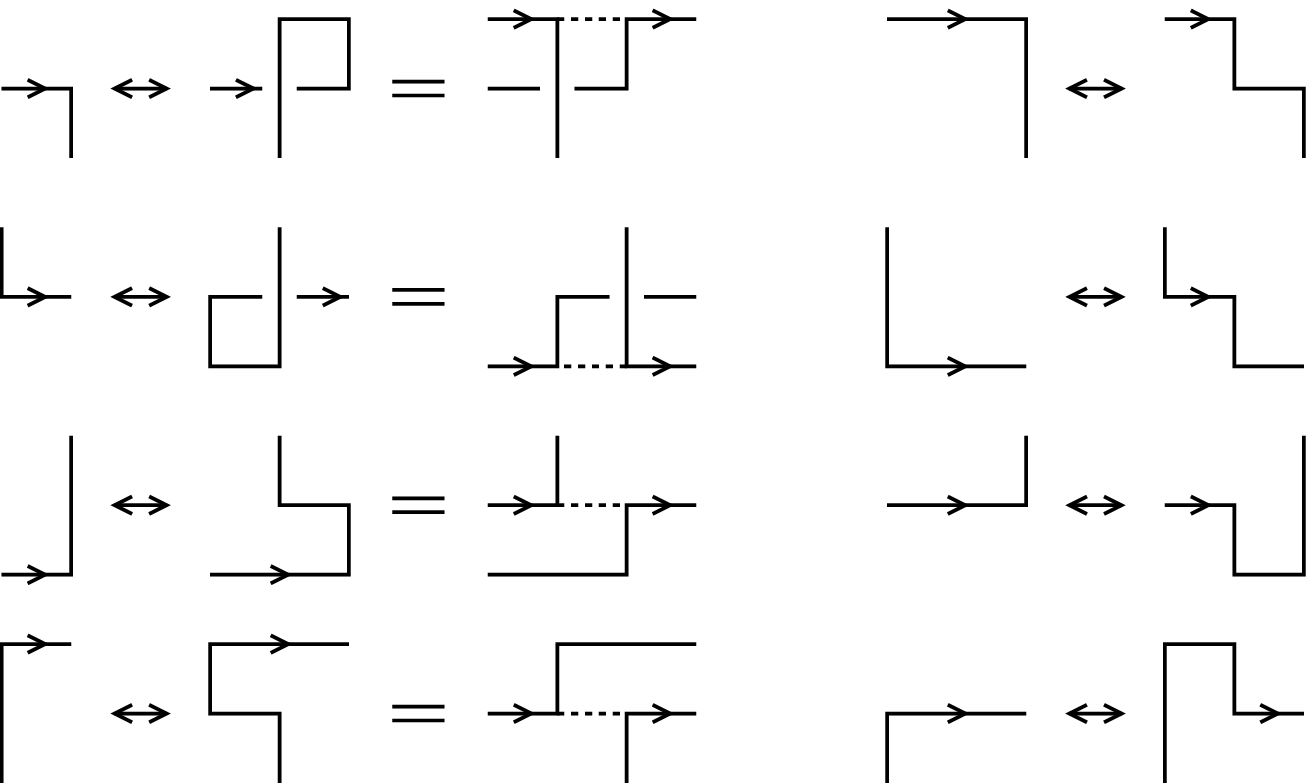}
\caption{} 
\label{legendrian}
\end{figure}

\begin{theorem} 
\label{legendrian-markov} 
{\rm (Markov Theorem for Legendrian links 
in $({\mathbb{R}}^3, \xi_{\rm sym})$)} \\ 
Let $D_1$ and $D_2$ be braided rectangular diagrams on $C_1$,  
and let $L_1$ and $L_2$ be Legendrian links corresponding to 
$D_1$ and $D_2$, respectively. 
Two Legendrian links $L_1$ and $L_2$ are Legendrian isotopic 
in $({\mathbb{R}}^3, \xi_{\rm sym})$ if and only if 
$D_1$ is obtained from $D_2$  
by a finite sequence of moves illustrated in Figure \ref{legendrian}. 
\end{theorem}

\begin{figure} 
 \includegraphics[keepaspectratio]{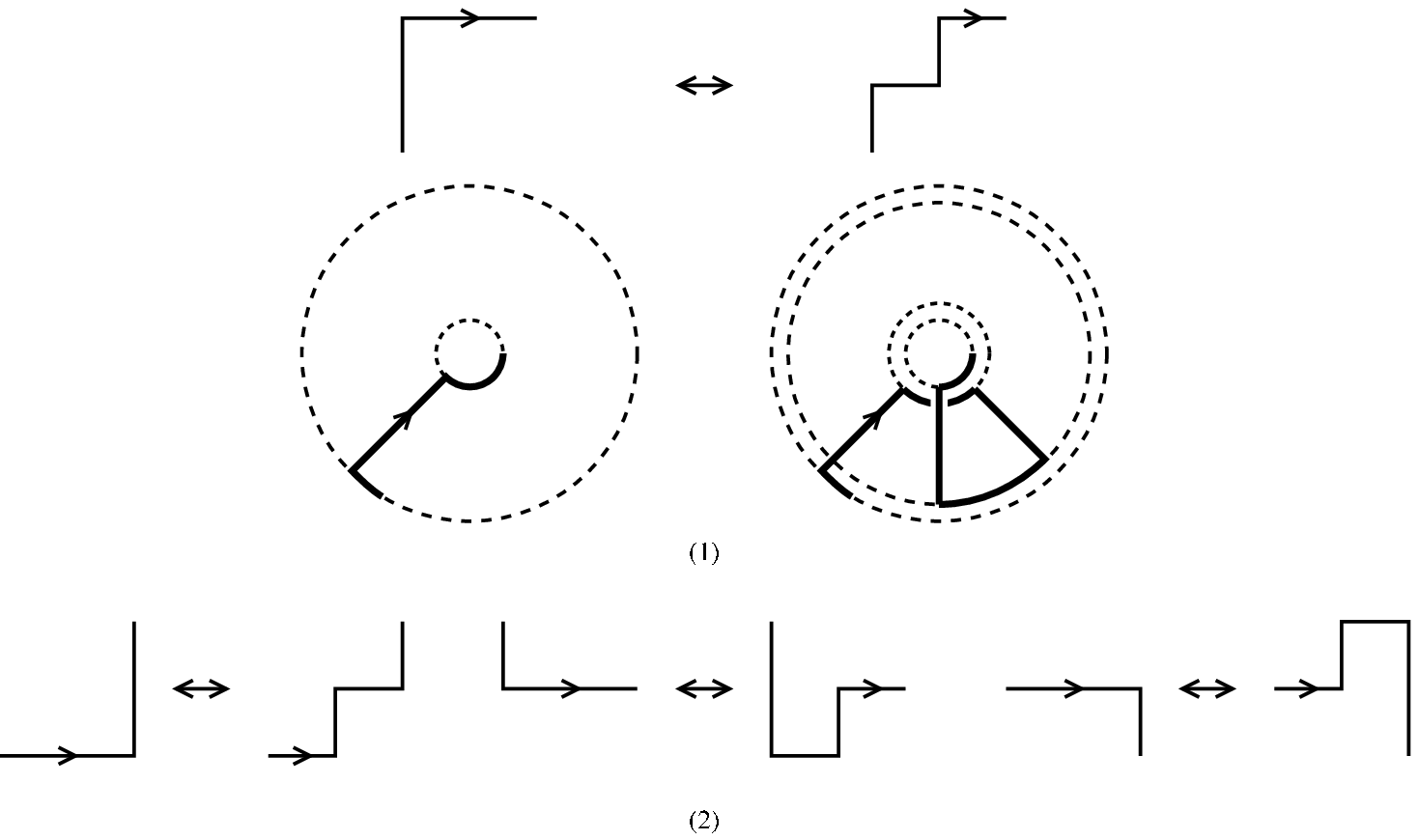}
\caption{} 
\label{transverse}
\end{figure}

\begin{remark} 
Let $p \colon {\mathbb{R}}^3 \to {\mathbb{R}}^2$ be 
a projection onto the $xy$-plane, and 
$L$ be a Legendrian link in $({\mathbb{R}}^3, \xi_{\rm sym})$. 
We denote $p(L)$ with over/under information by $p'(L)$. 
Let $D_1 \to D_2$ be a move in Figure \ref{transverse}. 
We note that $p'(L_1) \to p'(L_2)$ 
is a negative Reidemeister move of type I in knot theory, 
where $L_i$ $(i = 1, 2)$ is a Legendrian link corresponding 
to a rectangular diagram $D_i$. 
See Figure \ref{transverse} (1). 
This may be seen as a negative $``$local-(de)stabilization$"$. 
We note also that $T_+(L_1) \to T_+(L_2)$ 
is a transverse isotopy. 
See Theorem \ref{transverse-markov}. 
\end{remark}

\begin{figure} 
 \includegraphics[keepaspectratio]{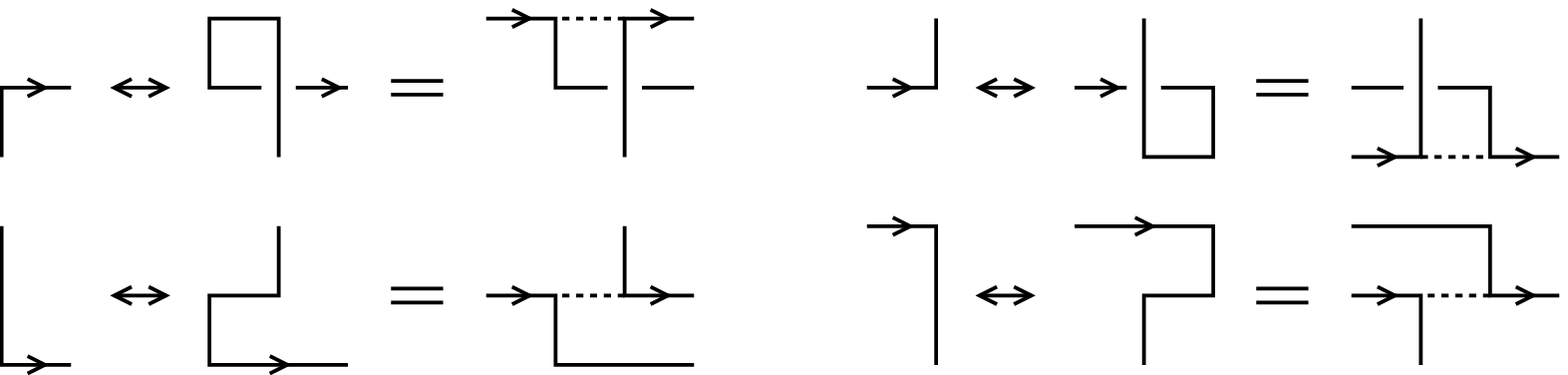}
\caption{} 
\label{transverse-not}
\end{figure}

\begin{remark} 
Let $D_1 \to D_2$ be a move in Figure \ref{transverse-not}, 
and $L_i$ $(i = 1, 2)$ be a Legendrian link 
corresponding to a rectangular diagram $D_i$. 
Then $T_+(L_1) \to T_+(L_2)$ is 
a negative (de)stabilization as closed braids, so
Theorem \ref{transverse-markov} shows that 
$T_+(L_1) \to T_+(L_2)$ is not a transverse isotopy. 
\end{remark}

Similar results as Theorems 
\ref{legendrian-reidemeister} and \ref{legendrian-markov} 
are obtained in \cite{ozsvath-szabo-thurston} 
in terms of rectangular diagrams, also known as grid-link diagrams.

\section{construction of Etnyre-Honda pair}

\begin{figure} 
 \includegraphics[keepaspectratio]{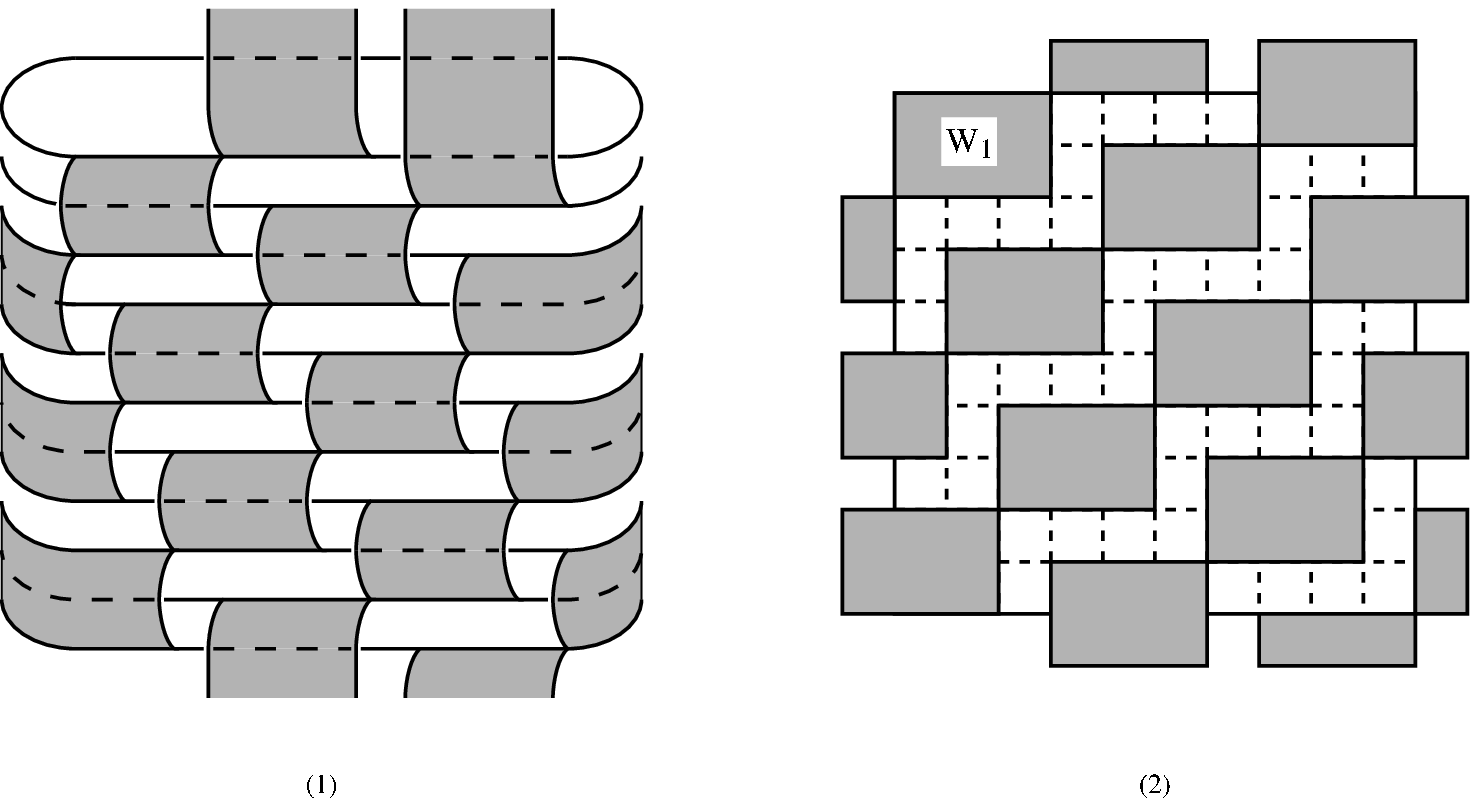}
\caption{} 
\label{torus}
\end{figure}

We start with a standardly embedded torus $U$ in ${\mathbb{R}}^3$. 
This torus $U$ may be described as 
$\{ (r, \theta, z) \ | \ (r-2)^2 + z^2 = 1 \}$ 
in the cylindrical coordinate of ${\mathbb{R}}^3$. 
Choose two sets of $pq+p+q$ numbers 
$z_1, \cdots, z_{pq+p+q}$ and 
$\theta_1, \cdots, \theta_{pq+p+q}$ 
satisfying the conditions 
$-1 < z_1 < \cdots < z_{pq+p+q} < 1$ and 
$0 < \theta_1 < \cdots < \theta_{pq+p+q} < \pi$. 
The intersection of $U$ with the plane $\{ (r, \theta, z) \ | \ z = z_i \}$ 
consists of two circles $L^i$ and $\ell^i$, 
where the $r$-coordinate of $L^i$ is larger than that of $\ell^i$. 
Let $M^i$ be the intersection of $U$ 
with the plane $\{ (r, \theta,z) \ | \ \theta = \theta_i \}$, which is a circle.  
The union of circles $\cup_{i=1}^{pq+p+q} (L^i \cup M^i)$ 
separate $U$ into $(pq+p+q)^2$ squares. 
Choose one rectangle $W_1$ on $U$ 
with $z_{pq+q} \leq z \leq z_{pq+p+q}$ and 
$\theta_1 \leq \theta \leq \theta_{q+1}$ occupying $pq$ squares. 
See Figure \ref{torus} (2). 
Choose a rectangle $W_{k+1}$ on $U$ 
occupying $pq$ squares in $p$ rows and $q$ columns 
so that the upper right corner of $W_k$ is identified with 
the lower left corner of $W_{k+1}$ for $k = 1, \cdots, pq+p+q$. 
Isotope $W_k$ so that 
$W_k \cap \{ (r, \theta, z) \ | \ (r-2)^2 + z^2 \leq 1 \}$ 
consists of a subarc of each of $L^i$ and $L^{i+p}$. 
See Figure \ref{torus} (1). 

\begin{figure} 
 \includegraphics[keepaspectratio]{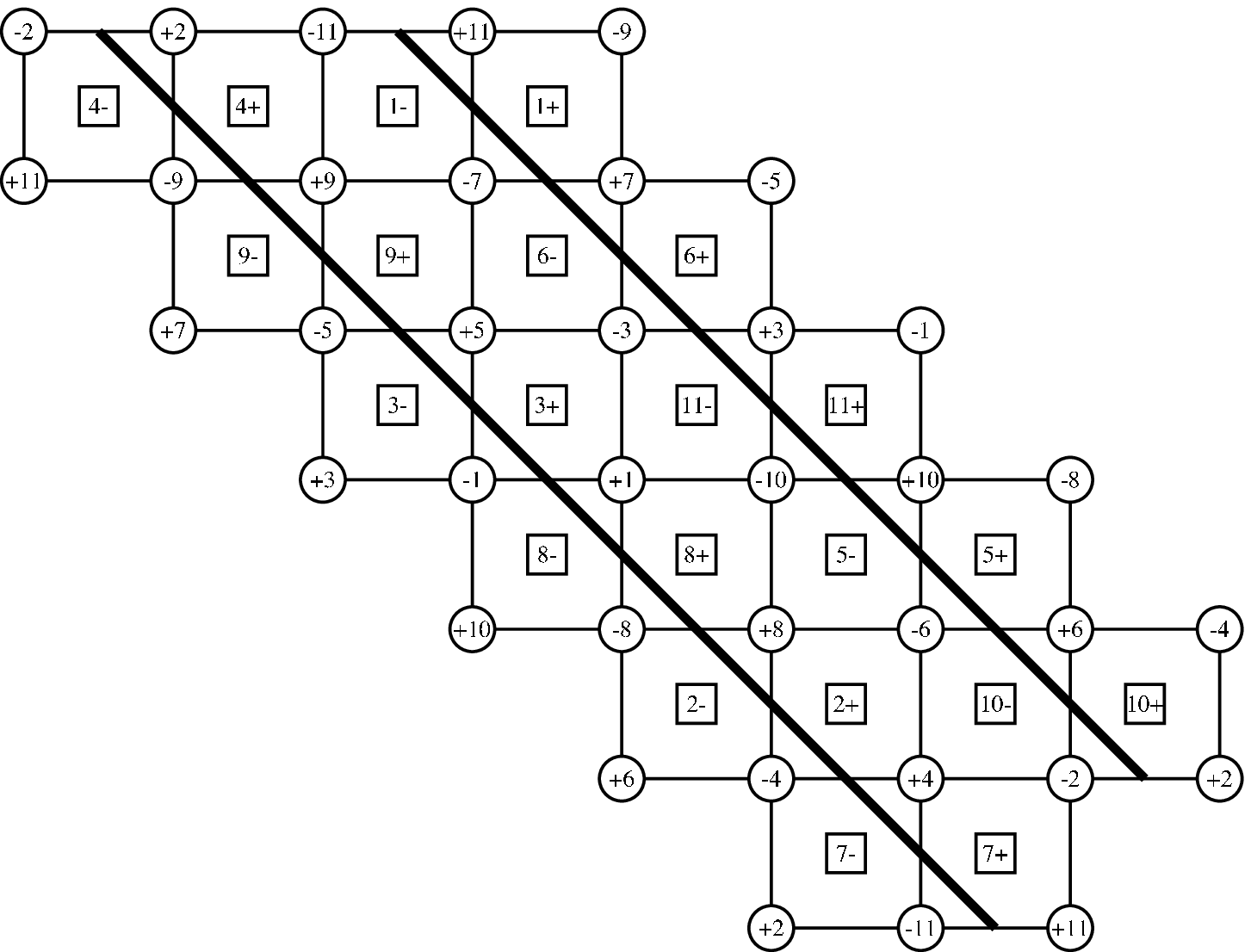}
\caption{} 
\label{tiling}
\end{figure}

Let $Z^i$ be the subdisc of the plane $\{ (r, \theta, z) \ | \ z = z_i \}$ 
with $\partial Z^i = L^i$. 
Let ${\mathcal{T}}$ denote a torus 
$\partial N_\varepsilon ((\cup_{i=1}^{pq+p+q} Z^i) \cup 
(\cup_{j=1}^{pq+p+q} W_j); S^3)$, 
that is the boundary of an $\varepsilon$-neighborhood of 
the branched surface 
$(\cup_{i=1}^{pq+p+q} Z^i) \cup (\cup_{j=1}^{pq+p+q} W_j)$ in $S^3$. 
Figure \ref{tiling} illustrates a tiling obtained from a braid foliation 
${\mathcal{T}} \cap \{ H(\theta) \}$ on ${\mathcal{T}}$, where 
the point with the mark $+k$ (resp. $-k$) represents 
the intersection of $\partial N(Z^k)$ and the $z$-axis 
with larger (resp. smaller) $z$-coordinate, and 
the point with the mark $k+$ (resp. $k-$) represents 
the hyperbolic singularity on the plane 
$\{ (r, \theta, z) \ | \ \theta = \theta_i + \varepsilon \}$ 
(resp. $\{ (r, \theta, z) \ | \ \theta = \theta_i - \varepsilon \}$) 
corresponding to $\partial N(W_j)$ (resp. $\partial N(W_{j+1})$). 
We notice that each of $G_{++}$ and $G_{--}$ consists of a circle, 
where $G_{++}$ and $G_{--}$ are the graphs 
defined in \cite{birman-finkelstein}.  
Let $C_1$ and $C_2$ be the annuli 
${\mathcal{T}} \setminus (G_{++} \cup G_{--})$, and 
let $c_1$ and $c_2$ be the core curve of $C_1$ and $C_2$, 
respectively. 
When $p=2$ and $q=3$, 
$c_1$ is a curve of slope $- \frac{2}{11}$ on ${\mathcal{T}}$ 
with respect to the coordinate system 
${\mathcal{C}}'$ in \cite{etnyre-honda}, 
where the boundary of the cabling annulus has slope $\frac{1}{0}$, 
and the meridian of ${\mathcal{T}}$ has slope $\frac{0}{1}$. 
Figure \ref{tiling} illustrates $c_1$ on ${\mathcal{T}}$. 

Next we look at ${\mathcal{T}}$ in $({\mathbb{R}}^3, \xi_{\rm sym})$. 
Isotope ${\mathcal{T}}$ so that 
$L^i$ has sufficiently small $r$ coordinate, and that 
the hyperbolic singularity on $\partial N(W_i)$ has 
sufficiently large $r$ coordinate. 
Each elliptic singularity in the characteristic foliation on ${\mathcal{T}}$ 
corresponds to an intersection of ${\mathcal{T}}$ with the $z$-axis, 
which is an elliptic singularity in a braid foliation. 
Each hyperbolic singularity in the characteristic foliation 
on ${\mathcal{T}}$ corresponds to 
a hyperbolic singularity on ${\mathcal{T}}$ in a braid foliation. 
Thus the characteristic foliation on ${\mathcal{T}}$ is 
isotopic to the corresponding braid foliation on ${\mathcal{T}}$. 

\begin{figure} 
 \includegraphics[keepaspectratio]{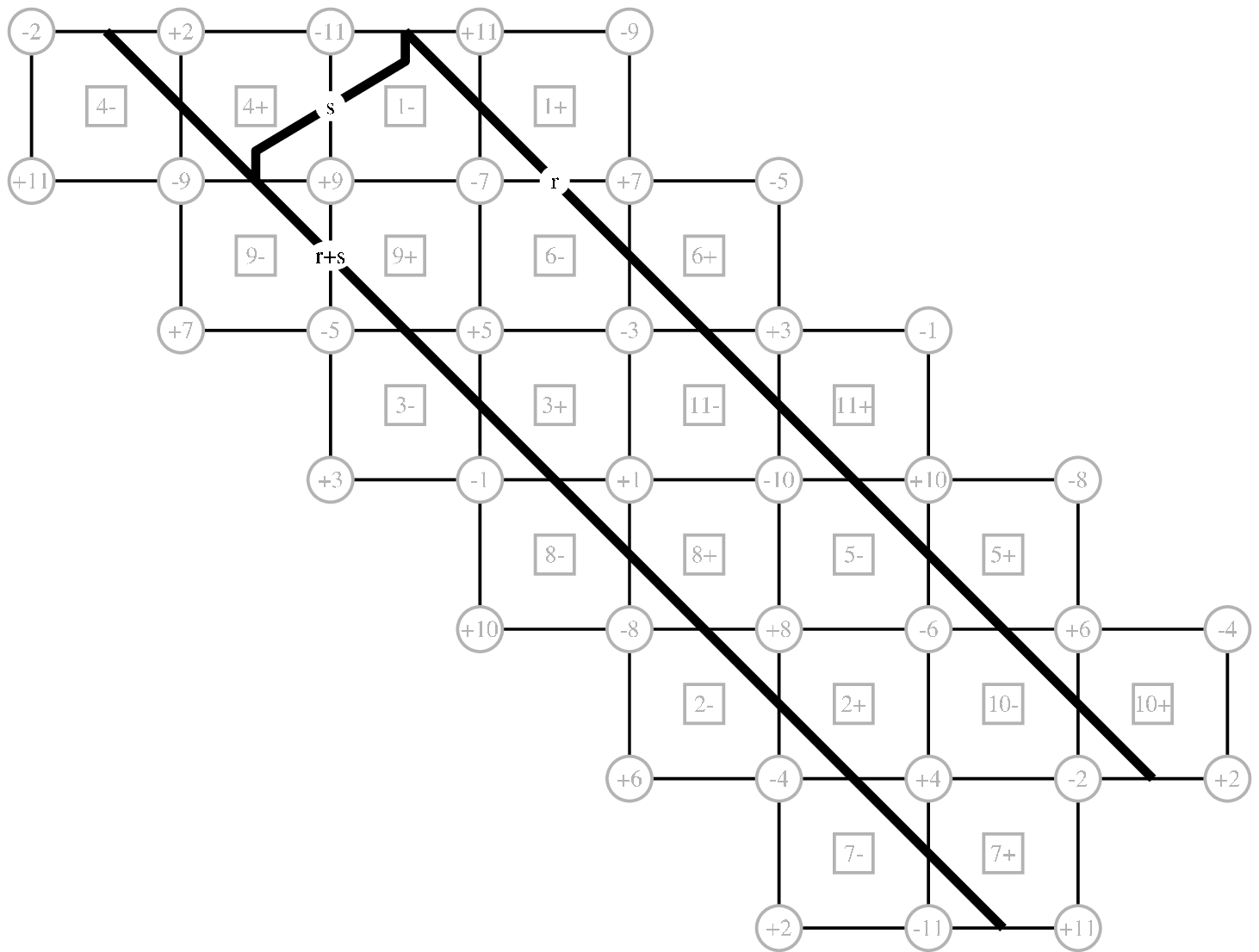}
\caption{} 
\label{rs-23}
\end{figure}

\begin{figure} 
 \includegraphics[keepaspectratio]{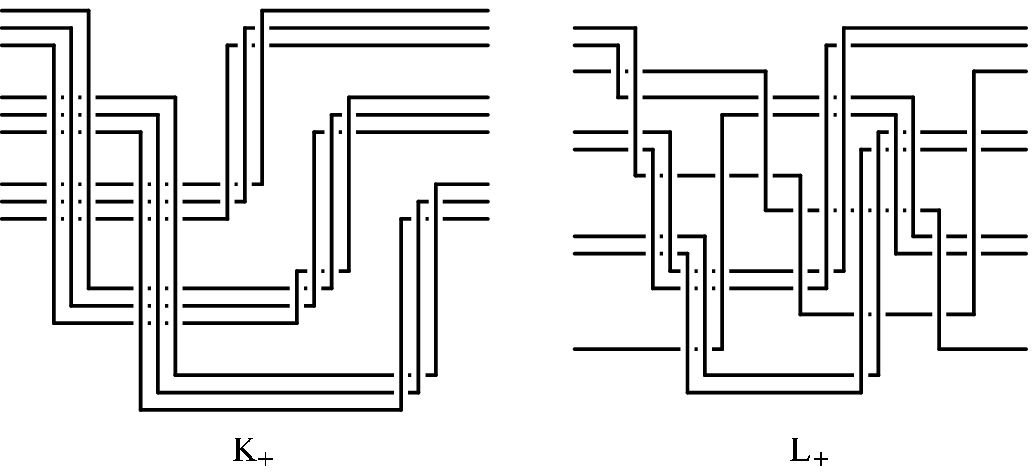}
\caption{} 
\label{pair}
\end{figure}

We may use the Giroux Elimination Lemma\cite{giroux} 
to isotope ${\mathcal{T}}$ 
in a small neighborhood of ${\mathcal{T}}$, and 
we may eliminate $G_{++}$ and $G_{--}$. 
Then ${\mathcal{T}}$ is a convex torus, and 
$c_1$ and $c_2$ are Legendrian divides on ${\mathcal{T}}$ 
with slope $-\frac{2}{11}$ with respect to ${\mathcal{C}}'$. 
Figure \ref{rs-23} illustrates a train-track on ${\mathcal{T}}$ 
constructed from $c_1$, $c_2$ and a small arc connecting them. 
Let $\ell(r, s)$ denote a simple closed curve supported by 
the train-track with weights $r$, $s$ and $r+s$, 
as illustrated in the figure. 
The topological knot type of $\ell(r, s)$ is 
a $(2r+s, r+s)$-cable of a $(2,3)$-torus knot. 
the Giroux Flexibility Theorem\cite{giroux} allows one to 
isotope ${\mathcal{T}}$ in a small neighborhood of ${\mathcal{T}}$ 
so that $\ell(r, s)$ is a Legendrian ruling curve on ${\mathcal{T}}$ 
with slope $-\frac{2r+s}{11r+5s}$ 
with respect to ${\mathcal{C}} '$ on ${\mathcal{T}}$. 
This Legendrian ruling curve 
may correspond to a braided rectangular diagram. 
When $r=s=1$, 
$\ell(1, 1)$ is $L_+$ in \cite{etnyre-honda}.  
A braided rectangular diagram of $\ell(1, 1)$ is obtained 
from the diagram in Figure \ref{pair}. 
It is easy to construct $K_+$ in \cite{etnyre-honda}, 
as illustrated in Figure \ref{pair}. 
A similar proof as in Lemma 6.3 in \cite{etnyre-honda} shows that 
$\ell(r, s)$ does not admit a Legendrian destabilization. 

\end{document}